\documentclass[11pt,reqno]{amsart}
\usepackage[T1]{fontenc}
\usepackage{lmodern}
\usepackage[margin=1in]{geometry}
\usepackage{mathtools,amssymb,amsthm}
\usepackage{microtype}
\usepackage{enumitem}
\usepackage{MnSymbol}
\usepackage{xcolor}
\usepackage[hidelinks]{hyperref}
\hypersetup{pdftitle={Counting fluctuations in the GOE-GUE crossover}}
\numberwithin{equation}{section}
\newtheorem{theorem}{Theorem}[section]
\newtheorem{proposition}[theorem]{Proposition}
\newtheorem{lemma}[theorem]{Lemma}
\newtheorem{corollary}[theorem]{Corollary}
\theoremstyle{remark}
\newtheorem{remark}[theorem]{Remark}
\newcommand{\E}{\mathbb{E}}
\newcommand{\Pp}{\mathbb{P}}
\newcommand{\R}{\mathbb{R}}
\newcommand{\C}{\mathbb{C}}
\newcommand{\scal}{{\operatorname{sc}}}
\newcommand{\dd}{\,\mathrm{d}}

\allowdisplaybreaks[1]
\title[Counting fluctuations in the GOE-GUE crossover]{Counting fluctuations in the GOE-GUE crossover}
\author{Connor Campbell}
\thanks{Department of Mathematics, University of Toronto,
Toronto, Ontario, Canada.}

\thanks{Email: \texttt{connork.campbell@mail.utoronto.ca}}

\address{Department of Mathematics, University of Toronto, Toronto, Ontario, Canada}

\email{connork.campbell@mail.utoronto.ca}

\date{\today}

\begin{document}
\begin{abstract}
We study bulk eigenvalue counting fluctuations in a Gaussian ensemble interpolating between the Gaussian orthogonal and Gaussian unitary ensembles. Realizing this ensemble as the distribution at a fixed time of Hermitian matrix Ornstein-Uhlenbeck flow started from the GOE, we prove a central limit theorem for the one-sided counting function, centered by the semicircle prediction. Its variance grows logarithmically, with a leading coefficient that interpolates between the GOE and GUE values. Specifically, at times $T_N \asymp N^{-1 + \vartheta}$ with fixed $0 < \vartheta < 1$, we prove that the one-sided bulk count is asymptotically Gaussian with leading variance $(1 - \vartheta/2)\pi^{-2}\log N$. Although the microscopic bulk correlations have their GUE limit throughout this regime, the leading counting variance retains a contribution from the initial GOE matrix.

%The proof follows the semicircle-centered logarithm of the characteristic polynomial along a deterministic characteristic. This approach separates the fluctuation inherited from the initial GOE matrix from an independent Gaussian fluctuation generated by the evolution, and identifies their respective contributions to the counting variance.
\end{abstract}
\maketitle

\section{Introduction and main result}\label{sec:Intro}
Let $\lambda_1 \leq \cdots \leq \lambda_N$ be the eigenvalues of a matrix from the Gaussian orthogonal ensemble (GOE) or the Gaussian unitary ensemble (GUE), normalized so that each off-diagonal entry has variance $N^{-1}$. For a bulk index $i = i_N$ with $i/N \to u \in (0, 1)$, let $\gamma_i$ be the semicircle quantile of order $i/N$. The central limit theorems of Gustavsson \cite{Gus} for the GUE and O'Rourke \cite{OR} for the GOE give $$\frac{N\rho_\scal(\gamma_i)(\lambda_i - \gamma_i)}{\sqrt{\log N}} \xrightarrow{\mathrm{d}} \operatorname{Normal}\left(0, \frac{1}{\beta\pi^2}\right),$$where $\rho_\scal$ is the semicircle density on $[-2, 2]$, and $\beta = 1$ for the GOE and $\beta = 2$ for the GUE. Thus, the fluctuation scale is the same in the two ensembles, but the limiting variances differ by a factor of two.

Pandey and Mehta \cite{PM} introduced the following Gaussian ensemble interpolating between the GOE and GUE. In our normalization, let $G_N$ and $V_N$ be independent $N \times N$ GOE and GUE matrices, and for $0\leq \alpha \leq 1$, define \begin{equation}\label{eq:ModelDefinition}
    W_N^{(\alpha)} := \sqrt{\frac{1 - \alpha^2}{1 + \alpha^2}}G_N + \sqrt{\frac{2\alpha^2}{1 + \alpha^2}}V_N.
\end{equation}
Then $W_N^{(0)} = G_N$ and $W_N^{(1)} = V_N$, and the off-diagonal variances remain equal to $N^{-1}$. Pandey and Mehta showed that the microscopic bulk correlations have their GOE limit when $\sqrt{N} \alpha_N \to 0$, and their GUE limit when $\sqrt{N}\alpha_N \to \infty$. A finite positive limit of $\sqrt{N}\alpha_N$ gives a nontrivial transition between the two. For small $\alpha$, the two coefficients in \eqref{eq:ModelDefinition} are $1 + O(\alpha^2)$ and $\sqrt{2}\alpha + O(\alpha^3)$, respectively. One may therefore think of $W_N^{(\alpha)}$ as a GOE matrix with a small independent GUE perturbation (specifically, $W_N \approx \mathrm{GOE} + \sqrt{2}\alpha \mathrm{GUE}$). Even when this perturbation tends to zero, the local correlations have their GUE limit provided $\sqrt{N}\alpha_N \to \infty$. This raises the question of whether individual eigenvalue fluctuations reach the GUE variance on the same scale as the local correlations. At first glance one might expect that the fluctuations of a single eigenvalue exhibit a similar transition from GOE to GUE fluctuations on the scale $\sqrt{N} \alpha_N \asymp 1$. However, this observable is not determined by the limiting microscopic correlation functions alone, and so such a result cannot be deduced from \cite{PM}. In fact, we will show that this naive guess is wrong, and that the leading counting variance interpolates between the GOE and GUE values across the scales $\alpha_N\asymp N^{-c}$, with fixed $0 < c < 1/2$, even though the microscopic correlations already have their GUE limit.

We study this question through the one-sided counting function $\mathcal{N}_N(E) := \#\{j : \lambda_j \leq E\}$, where the $\lambda_j$ now denote the eigenvalues of $W_N^{(\alpha_N)}$. The identity $$\{\lambda_i \leq E\} = \{\mathcal{N}_N(E) \geq i\}$$relates the fluctuations of an individual eigenvalue to those of the counting function. For the GOE and the GUE, the counting function at a fixed $E \in (-2, 2)$ satisfies a central limit theorem with leading variance $\pi^{-2}\log N$ and $(2\pi^2)^{-1}\log N$, respectively (see \cite{OR} and \cite{Gus}). We show that, in a regime where the local correlations already have their GUE limit, the coefficient of $\log N$ in the counting function variance can take any value strictly between the GOE and GUE coefficients.

\subsection{Main result} Writing $G_N = (g_{ij})$ and $V_N = (v_{ij})$, our normalization is $$g_{ij} \sim \operatorname{Normal}(0, N^{-1}) \quad (i < j), \qquad g_{ii} \sim \operatorname{Normal}(0, 2N^{-1}),$$and $$v_{ii} \sim \operatorname{Normal}(0, N^{-1}), \qquad \operatorname{Re}v_{ij}, \operatorname{Im}v_{ij} \sim \operatorname{Normal}(0, (2N)^{-1}) \quad (i < j).$$
All listed Gaussian coordinates are mutually independent, and the lower-triangular entries are determined by $g_{ji} = g_{ij}$ and $v_{ji} = \overline{v_{ij}}$. Write $$\rho_\scal(x) = \frac{1}{2\pi}\sqrt{(4 - x^2)_+}, \qquad F_\scal (E) = \int_{-\infty}^{E}\rho_\scal(x)\dd x $$and use the Stieltjes transform convention \begin{equation}\label{eq:semiCircleResolvent}
    m_\scal(z) = \int_\R \frac{\rho_\scal(x)}{x - z}\dd x = \frac{-z + \sqrt{z^2 - 4}}{2}, \quad z \in \C_+,
\end{equation}where the square root is chosen so that $\sqrt{z^2 - 4} \sim z$ at infinity. 
%For $0 \leq \alpha_N < 1$, it is convenient to parametrize the model by $$T_N = \log\frac{1 + \alpha_N^2}{1 - \alpha_N^2}, \qquad \alpha_N^2 = \tanh(T_N/2).$$As will be shown in Section~\ref{sec:Characteristic}, this is the time at which a Hermitian matrix Ornstein-Uhlenbeck flow started from $G_N$ has the law of $W_N^{(\alpha_N)}$. 
Our main result is the following:
\begin{theorem}
    \label{thm:MainResult}Fix $E_0 \in (-2, 2)$ and $0 < \vartheta < 1$. Let $T_N$ be deterministic with $T_N \asymp N^{-1 + \vartheta}$, and set $\alpha_N^2 = \tanh(T_N/2)$. Let $\mathcal{N}_N(E_0)$ count the eigenvalues of $W_N^{(\alpha_N)}$ in $(-\infty, E_0]$. Then \begin{equation}\label{eq:MainResult}
    \frac{\mathcal{N}_N(E_0) - NF_\scal(E_0)}{\pi^{-1}\sqrt{(1 - \vartheta/2)\log N}} \xrightarrow{\mathrm{d}} \operatorname{Normal}(0, 1). 
\end{equation}Moreover, \begin{align}
    \E\mathcal{N}_N(E_0) &= NF_\scal(E_0) + O((\log N)^{1/4}), \label{eq:MeanAsymptotic} \\
    \operatorname{Var}(\mathcal{N}_N(E_0)) &= \frac{1 - \vartheta/2}{\pi^2}\log N + O\left((\log N)^{3/4}\right). \label{eq:VarianceAsymptotic}
\end{align}In particular, the central limit theorem also holds after centering by $\E\mathcal{N}_N(E_0)$ and dividing by the exact standard deviation. 
\end{theorem}

In this regime, $\sqrt{N}\alpha_N \asymp N^{\vartheta/2} \to \infty$, so the local correlations have their GUE limit for every fixed $\vartheta \in (0, 1)$. The coefficient in \eqref{eq:VarianceAsymptotic}, however, remains strictly between the GOE and GUE values.

\begin{remark}[Multipoint fluctuations] \label{rem:MultipointFluctuations} In this paper, we restrict attention to the counting function at a single energy. However, the same methods extend to joint fluctuations at any fixed number of bulk energies. For example, one may consider two deterministic energies, in a fixed compact subinterval of $(-2, 2)$, satisfying $|E_1 - E_2| = (1 + o(1))N^{-a}$ for $0 < a < 1$. 
    The additional calculation needed to identify the joint Gaussian limit is to compute the covariance of the two counts as a function of $a$ and $\vartheta$. We outline how the proof extends and give the leading covariance in Section~\ref{sec:RemoveRegularization}.
\end{remark}

\subsection{Related work}Bourgade and Mody \cite[Theorem 1.6 and Proposition B.1]{BM} extended the individual eigenvalue and counting central limit theorems from the GOE and GUE to real symmetric and complex Hermitian Wigner matrices under a subgaussian tail assumption. The limiting variances agree with those of the GOE and GUE in the real and complex cases, respectively. In the model considered here, the relative variances of the real and imaginary parts of the entries depend on $N$, and the leading counting variance can lie strictly between the GOE and GUE values. 

For Gaussian models interpolating between symmetry classes, Pandey and Mehta \cite{PM} computed the bulk correlation functions for the model \eqref{eq:ModelDefinition} considered in this paper. Forrester, Nagao, and Honner \cite{ForresterNagaoHonner} computed limiting correlation functions at the soft edge for Gaussian ensembles interpolating between the GOE and GUE. They also obtained hard-edge limits for ensembles interpolating between the Laguerre orthogonal and unitary ensembles. 

He \cite{He} proved a central limit theorem for smooth mesoscopic linear statistics of Wigner matrices, including the GOE-GUE interpolation considered here. He showed that the transition between the GOE and GUE limiting variances occurs when the evolution time is comparable to the mesoscopic scale of the test function. His result does not cover the one-sided counting function considered here, whose test function is discontinuous.

\subsection{Proof strategy}
Write $T=T_N$, where $T_N$ is defined as in Theorem~\ref{thm:MainResult}. We study $W_N^{(\alpha_N)}$ through a matrix process $H_t$ that starts from $G_N$ and has the law of $W_N^{(\alpha_N)}$ at time $T$. This process is driven by Hermitian Brownian motion independent of $G_N$ and is defined in Section~\ref{sec:Characteristic}. Denote its eigenvalues by $\lambda_i(t)$.

The main quantity we study in the proof is the centered logarithm of the characteristic polynomial, $$L_t(z) := \sum_{i=1}^N\operatorname{Log}_-(\lambda_i(t) - z) -N\int_\R\operatorname{Log}_-(x - z)\rho_\scal(x)\dd x, \qquad z\in\C_+,$$
where $\operatorname{Log}_-$ has argument in $(-\pi,0)$. For $\eta >0$,  one sees that $-\pi^{-1}\operatorname{Im}\operatorname{Log}_-(x-E-\mathrm{i}\eta)$ is a smooth approximation to $\mathbf{1}_{\{x\leq E\}}$, and so the quantity $-\pi^{-1}\operatorname{Im}L_T(E+\mathrm{i}\eta)$
approximates the centered counting function $\mathcal{N}_N(E) - NF_\scal(E)$. We take
$\eta_*=(\log N)^{1/4}/N$ and prove in Lemma~\ref{lem:DesmoothingEstimate} that $$\left\|\mathcal N_N(E_0)-NF_\scal(E_0) +\frac{1} {\pi}\operatorname{Im}L_T(E_0+\mathrm{i}\eta_*)\right\|_{L^2} = O((\log N)^{1/4}).$$ It therefore suffices to study $\operatorname{Im}L_T(E_0+\mathrm{i}\eta_*)$ up to errors of size $o ( \sqrt{ \log N})$.

We will study $L_T$ via the method of characteristics (see, e.g., \cite{HuangLandon}). We will promote $z$ to a time dependent quantity satisfying the equations $$\dot z_t = -m_\scal(z_t) - \frac{z_t}{2}, \qquad z_T = E_0 + \mathrm{i}\eta_*.$$ Analyzing the evolution of $L_t$ along this characteristic, we obtain $$L_T(z_T) = L_0(z_0) + M_N -\frac{N}{2}\int_0^T \bigl(m_t(z_t) - m_\scal(z_t)\bigr)^2\dd t,$$where $m_t(z) = N^{-1}\sum_i(\lambda_i(t)-z)^{-1}$, and $M_N$ is a stochastic integral with integrands $N^{-1/2}(\lambda_i(t) - z_t)^{-1}$. 
%We will show that the drift term, $\frac{N}{2}\int_0^T\bigl(m_t(z_t) - m_\scal(z_t)\bigr)^2\dd t$, is $o(\sqrt{\log N})$ in $L^2$, and that $\operatorname{Im}M_N$ can be approximated by a Gaussian variable independent of $G_N$.}

We next divide $[0, T]$ into two intervals, $[0, t_1]$ and $[t_1, T]$, where $t_1$ is chosen so that $\eta_* \ll \operatorname{Im}z_{t_1} \ll T$. This choice of $t_1$ serves two purposes. First, for $t \in [0, t_1]$, $\operatorname{Im}z_t$ is sufficiently large for the eigenvalue rigidity estimates of \cite{EYY} to show that the drift term on this interval, $\frac{N}{2}\int_0^{t_1}\bigl(m_t(z_t) - m_\scal(z_t)\bigr)^2\dd t$, is $o(\sqrt{\log N})$ in $L^2$. Second, these estimates allow us to replace $\lambda_i(t)$ by the deterministic semicircle quantiles $\gamma_i$ in the integrand defining $M_N$, with an error of $o(1)$ in $L^2$. After this replacement, the integrands are deterministic, so the imaginary part of the resulting integral is Gaussian and independent of $G_N$. Its variance can then be computed explicitly. 

In contrast, on $[t_1, T]$, the imaginary part of $z_t$ decreases to $\eta_*$, where the rigidity estimates alone are no longer sufficient to control the two terms. In this regime, we instead consider the evolution equations for $m_t - m_\scal$ and its derivative. Starting from the estimates at $t_1$, we control the martingale increments and use a continuity argument to obtain stronger bounds which hold for heights down to $\operatorname{Im}z_T = \eta_*$. These bounds show that the contributions of $t \in [t_1, T]$ to both the drift and $M_N$ terms are negligible on the $\sqrt{\log N}$ scale.

Denote the Gaussian variable obtained from the integral over $[0, t_1]$ by $Z_N$, and set $X_N := \operatorname{Im}L_0(z_0)$. The previous estimates give \begin{equation}\label{eq:IntroFluctuationDecomposition}
    \operatorname{Im}L_T(z_T)= X_N + Z_N + R_N,  \qquad \|R_N\|_{L^2} = o(\sqrt{\log N}).
\end{equation}The initial term $X_N$ is evaluated at height $\operatorname{Im}z_0 \asymp T \asymp N^{-1 + \vartheta}$, and is therefore a mesoscopic linear statistic of the GOE. Its Gaussian limit  follows from standard methods for mesoscopic linear statistics (here we apply the estimate in \cite[Appendix B]{BLZ}). Together with a direct computation of the variance of $Z_N$, and using the independence of $X_N$ and $Z_N$, we obtain $$\frac{X_N + Z_N}{\sqrt{\log N}} \xrightarrow{\dd}\operatorname{Normal}(0, 1 - \vartheta/2).$$
Since the remainder in \eqref{eq:IntroFluctuationDecomposition} and the error in approximating the counting function are both $o(\sqrt{\log N})$ in $L^2$, we obtain the central limit theorem and the leading variance in Theorem~\ref{thm:MainResult}.

\subsection{Notation and organization.} Constants $C, c > 0$ are independent of $N$ and may change from line to line. They may depend on $E_0, \vartheta,$ the implicit constants in $T_N \asymp N^{-1 + \vartheta}$ and the fixed choice of $A$ introduced in Section~\ref{sec:Rigidity}. All asymptotic assertions are as $N \to \infty$. We use $f_N \asymp g_N$ to mean $c_*f_N \leq g_N \leq C_*f_N$ for $N$ large, and $f_N \ll g_N$ to mean $f_N/g_N \to 0$. For a complex random variable $X$, we write $\|X\|_{L^2} = \left(\E |X|^2\right)^{1/2}$. For a real continuous local martingale $Q$, $\langle Q \rangle$ denotes its quadratic variation. We use $\|\cdot\|_{\text{op}}$ for the operator norm and the convention $\inf\emptyset = \infty$. An event holds with overwhelming probability when its complement has probability at most $N^{-D}$ for every fixed $D > 0$ and all sufficiently large $N$. 

Section~\ref{sec:Characteristic} gives the characteristic and exact evolution identities. Sections~\ref{sec:Rigidity}-\ref{sec:SecondMomentBounds} establish the coarse and fine resolvent estimates. Section~\ref{sec:GaussianSection} constructs the independent Gaussian approximation and computes its variance. Section~\ref{sec:InitialGOEField} treats the initial GOE field, and Section~\ref{sec:RemoveRegularization} removes the regularization and proves Theorem~\ref{thm:MainResult}.  

\section{Characteristic and exact evolution}\label{sec:Characteristic}
Let $H_t$ be the solution to \begin{equation}\label{eq:OrnsteinUhlenbeck}
    \dd H_t = N^{-1/2}\dd B_t - \frac{1}{2}H_t\dd t, \qquad H_0 = G_N,
\end{equation}where $B_t$ is a complex Hermitian matrix Brownian motion independent of $G_N$. Specifically, we use the normalization $$(B_t)_{ii} = \beta_{ii}(t), \qquad (B_t)_{ij} = \frac{\beta_{ij}(t) + \mathrm{i}\widetilde{\beta}_{ij}(t)}{\sqrt{2}} \quad (i < j),$$where the $\beta_{ii}(t)$, $\beta_{ij}(t)$, and $\widetilde{\beta}_{ij}(t)$ are independent standard real Brownian motions. Let $(\mathcal{F}_t, t \geq 0)$ denote the usual augmentation of the filtration generated by $H_0$ and $B$. By the method of integrating factors, \eqref{eq:OrnsteinUhlenbeck} has an explicit solution equal in distribution to \begin{equation}\label{eq:DynamicalModel}
    H_t \overset{\dd}{=} e^{-t/2}G_N + \sqrt{1 - e^{-t}}V_N
\end{equation}for each fixed $t \geq 0$, where $V_N$ is a GUE matrix independent of $G_N$. Recalling $T = T_N$ defined in Theorem~\ref{thm:MainResult}, we see that $H_T$ has the same distribution as $W_N^{(\alpha_N)}$. We therefore work with $H_T$ throughout the proof. 

Denote the ordered eigenvalues by $\lambda_1(t) < \cdots < \lambda_N(t)$, and write $\mathcal{N}_N(E;t) := \#\{i : \lambda_i(t) \leq E\}$ so that $\mathcal{N}_N(E) = \mathcal{N}_N(E;T)$. For $z \in \C_+$, set  \begin{gather}
    m_t(z) := \frac{1}{N}\sum_{i = 1}^N\frac{1}{\lambda_i(t) - z}, \qquad  \delta_t(z) :=m_t(z) - m_\scal(z), \label{eq:resolventDefinition} \\
    L_t(z) := \sum_{i = 1}^N\operatorname{Log}_-(\lambda_i(t) - z) - N\int_\R \operatorname{Log}_-(x-z)\rho_\scal(x)\dd x .\label{eq:LogFieldDefinition}
\end{gather}
Here, $\operatorname{Log}_-$ is the branch of the logarithm on the lower half-plane with argument in $(-\pi, 0)$. Differentiating \eqref{eq:LogFieldDefinition} gives
\begin{equation}\label{eq:LogFieldzDerivative}
    \partial_zL_t(z) = -N\delta_t(z). 
\end{equation}

Set $$\eta_* := \frac{(\log N)^{1/4}}{N}, \qquad b(z) := m_\scal(z) + \frac{z}{2},$$and let $z_t = x_t + \mathrm{i}y_t$ solve \begin{equation}\label{eq:characteristicDefinition}
    \dot{z}_t = -b(z_t), \quad z_T = E_0 + \mathrm{i}\eta_*.
\end{equation}We will use the following estimates for this characteristic. 
\begin{lemma}\label{lem:CharacteristicGeometry} For all sufficiently large $N$, the characteristic lies in $\C_+$ for $0 \leq t \leq T$, and its real part lies in a fixed compact subinterval of $(-2, 2)$. Uniformly for $0\leq t \leq T$, \begin{align}
    c \leq -\dot{y_t} = \operatorname{Im}b(z_t) &\leq C, & \eta_* + c(T - t) \leq y_t &\leq \eta_* + C(T - t), \label{eq:boundsImaginaryCharacteristic} \\
    |z_t - z_T| &\leq C(T - t), & y_0 &\asymp T. \label{eq:CharacteristicDisplacement}
\end{align}
For $0 \leq r \leq s \leq T$ and $a > 1$, \begin{equation}\label{eq:characteristicIntIneq}
    \int_r^sy_u^{-a}\dd u \leq C_a y_s^{1 - a}, \qquad \int_r^s \frac{\dd u}{y_u} \leq C\log\frac{y_r}{y_s}. 
\end{equation}
Moreover, \begin{equation}\label{eq:characteristicIntEq}
    \int_r^s\frac{\operatorname{Im}b(z_u)}{y_u}\dd u = \log\frac{y_r}{y_s}.
\end{equation}At time zero, \begin{equation}\label{eq:logInitialHeight}
    x_0 = E_0 + O(T), \quad \log\frac{1}{y_0} = (1 - \vartheta)\log N + O(1). 
\end{equation}\end{lemma}
\begin{proof} 
    By \cite[Equation (1.11)]{BourgadeExtremeGaps}, the solution of \eqref{eq:characteristicDefinition} is $$z_t = z_T\cosh\frac{T - t}{2} + \sqrt{z_T^2 - 4}\sinh\frac{T - t}{2}, \qquad 0 \leq t \leq T,$$where the square root has positive imaginary part. Since $E_0 \in (-2, 2), T \to 0$, and $\eta_*/T \to 0$, this formula gives \eqref{eq:boundsImaginaryCharacteristic}-\eqref{eq:CharacteristicDisplacement}. Taking $t = 0$ and using $T \asymp N^{-1 + \vartheta}$ gives \eqref{eq:logInitialHeight}. The change of variables $y = y_u$, together with $-\dot{y}_u \geq c$, proves \eqref{eq:characteristicIntIneq}. Integrating $-\dot{y}_u /y_u = \operatorname{Im}b(z_u)/y_u$ gives \eqref{eq:characteristicIntEq}. Finally, we also note that \eqref{eq:semiCircleResolvent} shows that each fixed order derivative of $m_\scal$ is uniformly bounded along the characteristic. 
\end{proof}
We next record the evolution of the Stieltjes transform and the logarithmic field along $z_t$. Write $$\delta_t := \delta_t(z_t), \qquad q_t := m'_t (z_t) - m'_{\mathrm{sc}} (z_t). % \partial_z\delta_t(z_t).
$$
%{\cor I wrote out $\partial_z \delta$ more explicitly, becuase if you write $\partial_z f (z_t)$ then it could be interpreted as the $\partial_z$ derivative also hitting $z_t$ if you think of the characteristic map $z_t$ as a function of $z$}
\begin{lemma}\label{lem:EigenvalueDBMandSDEs}
There exists a standard $N$-dimensional Brownian motion $(B_1, ..., B_N)$ with respect to $(\mathcal{F}_t)$, independent of $H_0$, such that \begin{align}
    \dd \delta_t &= \left[(m_\scal'(z_t) + \tfrac{1}{2})\delta_t + \delta_tq_t\right]\dd t + \dd M_t^{(0)}, \label{eq:deltaSDE} \\
    \dd q_t &= \left[q_t^2 + (2m_\scal'(z_t) + 1)q_t + \delta_t\partial_z^2m_t(z_t)\right]\dd t + \dd M_t^{(1)}, \label{eq:qSDE}
\end{align}where $M_0^{(0)} = M_0^{(1)} = 0$ and \begin{equation}
    \dd M_t^{(0)} = -\frac{1}{N^{3/2}}\sum_{i}\frac{\dd B_i(t)}{(\lambda_i(t) - z_t)^2}, \qquad \dd M_t^{(1)} = -\frac{2}{N^{3/2}}\sum_i\frac{\dd B_i(t)}{(\lambda_i(t) - z_t)^3}. \label{eq:MartingalesInSDEs}
\end{equation}Moreover, \begin{equation}
    L_T(z_T) = L_0(z_0) + M_N - \frac{N}{2}\int_0^T\delta_t^2\dd t, \qquad M_N = \frac{1}{\sqrt{N}}\sum_i\int_0^T\frac{\dd B_i(t)}{\lambda_i(t) - z_t}. \label{eq:LogFieldTerminalTime}
\end{equation}

%\begin{equation}\label{eq:EigenvalueDBMSDE} \dd x_i(t) = \frac{1}{\sqrt{N}}\dd B_i(t) + \left(\frac{1}{N}\sum_{j \neq i}\frac{1}{x_i(t) - x_j(t)} - \frac{x_i(t)}{2}\right)\dd t. \end{equation} Set $\delta_t = \delta_t(z_t)$ and $q_t = \partial_z\delta_t(z_t)$. Then, \begin{align} \dd \delta_t &= \left[(m_\scal'(z_t) + \tfrac{1}{2})\delta_t + \delta_tq_t\right]\dd t + \dd M_t^{(0)}, \label{eq:deltaSDE} \\ \dd q_t &= \left[q_t^2 + (2m_\scal'(z_t) + 1)q_t + \delta_t\partial_z^2m_t(z_t)\right]\dd t + \dd M_t^{(1)}, \label{eq:qSDE} \end{align}where $M_0^{(0)} = M_0^{(1)} = 0$ and \begin{equation}\dd M_t^{(0)} = -\frac{1}{N^{3/2}}\sum_{i}\frac{\dd B_i(t)}{(x_i(t) - z_t)^2}, \qquad \dd M_t^{(1)} = -\frac{2}{N^{3/2}}\sum_i\frac{\dd B_i(t)}{(x_i(t) - z_t)^3}. \label{eq:MartingalesInSDEs}\end{equation}Finally, \begin{equation} L_T(z_T) = L_0(z_0) + M_N - \frac{N}{2}\int_0^T\delta_t^2\dd t, \qquad M_N = \frac{1}{\sqrt{N}}\sum_i\int_0^T\frac{\dd B_i(t)}{x_i(t) - z_t}. \label{eq:LogFieldTerminalTime}\end{equation}
\end{lemma}
\begin{proof} This is a direct calculation. For completeness, we indicate the relevant references. The Brownian motions are obtained from the eigenvalue process as in \cite[Section 4.3.1]{AGZ}. Their independence of $H_0$ follows since they form an $(\mathcal{F}_t)$-Brownian motion and $H_0$ is $\mathcal{F}_0$-measurable. 

Applying \cite[Equation (3.19)]{HuangLandon} with $\beta = 2$ and $V(x) = x^2/2$, and using $$bm_\scal' + \frac{1}{2}m_\scal = 0,$$gives \eqref{eq:deltaSDE} along the characteristic \eqref{eq:characteristicDefinition}. Differentiating the evolution equation for $m_t(z)$, $$\dd m_t(z) = \left[(m_t + \tfrac{z}{2})\partial_zm_t + \tfrac{1}{2}m_t\right]\dd t - \frac{1}{N^{3/2}}\sum_i\frac{\dd B_i(t)}{(\lambda_i(t) - z)^2},$$ with respect to $z$, then evaluating along $z_t$ and using $\dot{z}_t = -b(z_t)$ gives the evolution of $q_t$. The identity $$bm_\scal'' + (m_\scal')^2 + m_\scal' = 0,$$obtained by differentiating $bm_\scal' + m_\scal/2 = 0$, yields \eqref{eq:qSDE}. Finally, \eqref{eq:LogFieldTerminalTime} follows by setting $\beta = 2$ in \cite[Equations (7.13)-(7.17)]{LandonSosoe}, multiplying by $N$ to match our definition of $L_t$, and integrating over $[0, T]$.  
\end{proof}
We also make use of the following standard deterministic inequalities. For any Stieltjes transform $$m(z) := \frac{1}{N}\sum_{i = 1}^N\frac{1}{\lambda_i - z}, \qquad z \in \mathbb{C}_+, \quad \lambda_i \in \R,$$writing $z = E + \mathrm{i}\eta$ and using $|\lambda_i - z| \geq \eta$, we obtain \begin{equation}\label{eq:DeterministicResolventInequalities}
    \operatorname{Im}m(z) = \frac{1}{N}\sum_i\frac{\eta}{|\lambda_i - z|^2}, \qquad |m'(z)|\leq \frac{\operatorname{Im}m(z)}{\eta}, \qquad |m''(z)| \leq \frac{2\operatorname{Im}m(z)}{\eta^2}. 
\end{equation}Applying these estimates to \eqref{eq:MartingalesInSDEs} yields \begin{align}
    \dd \langle \operatorname{Re}M^{(0)}\rangle_t, \ \dd \langle \operatorname{Im} M^{(0)}\rangle_t &\leq \frac{\operatorname{Im}m_t(z_t)}{N^2y_t^3}\dd t \label{eq:QuadVariationBoundDeltaMG} \\
    \dd \langle \operatorname{Re}M^{(1)}\rangle_t, \ \dd \langle \operatorname{Im} M^{(1)}\rangle_t &\leq \frac{4\operatorname{Im}m_t(z_t)}{N^2y_t^5}\dd t. \label{eq:QuadVariationBoundqMG}
\end{align}

\section{Uniform rigidity estimates}\label{sec:Rigidity} 
Let $\gamma_i$ denote the semicircle quantiles, defined by $F_\scal(\gamma_i) = i/N$ for $1 \leq i < N$, and $\gamma_N = 2$. Fix $A > 0$ sufficiently large, as required in Proposition~\ref{prop:rigidityEvent}, and set \begin{equation}\label{eq:intermediateHeightParams}
    r_N := \exp(A(\log \log N)^2), \qquad h_N := \frac{r_N^4}{N}. 
\end{equation}It follows from \eqref{eq:intermediateHeightParams} and Lemma~\ref{lem:CharacteristicGeometry} that, for all sufficiently large $N$, there is a unique deterministic $t_1 \in (0, T)$ such that $y_{t_1} = h_N$. Moreover, we have \begin{equation}\label{eq:hNScales}
    \eta_* \ll h_N \ll y_0, \qquad T - t_1 \leq Ch_N, \qquad \log(h_N/\eta_*) = O(\log r_N) = o(\log N). 
\end{equation}

This choice of $h_N$ is motivated by the approximation of the martingale in \eqref{eq:LogFieldTerminalTime}. On $[0, t_1]$, we will replace $\lambda_i(t)$ by $\gamma_i$ in its integrand. For indices with $\gamma_i$ in a fixed compact subinterval of $(-2, 2)$, Proposition~\ref{prop:rigidityEvent} will give $|\lambda_i(t) - \gamma_i|\leq Cr_N/N$ uniformly in $t$ on an event with overwhelming probability. Since $y_t \geq h_N$ for $t \in [0, t_1]$, we have \begin{equation}\label{eq:RigidityVsHeightRatio}
    \frac{r_N/N}{h_N} = \frac{r_N}{Nh_N} = \exp(-3A(\log \log N)^2) = o(1). 
\end{equation}This allows us to estimate the difference between $(\lambda_i(t) - z_t)^{-1}$ and $(\gamma_i - z_t)^{-1}$. The remaining martingale integral over $[t_1, T]$ will be estimated separately, using the bounds in \eqref{eq:hNScales}. We first record the martingale inequalities used below. \vspace{2mm}\newline
{\it Martingale estimates.} We use the upper Burkholder-Davis-Gundy inequality with explicit dependence on the moment. For a real continuous local martingale $Q$ on $[a, b]$ with $Q_a = 0$, write$$Q_b^* := \sup_{a \leq t \leq b}|Q_t|, \qquad [Q]_{a, b} := \langle Q \rangle_b - \langle Q \rangle_a.$$
By \cite[Theorem 5.16 and its proof]{LeGall}, for every real $p \geq 2$,  \begin{equation}\label{eq:LpBDGBound}
    \|Q_b^*\|_{L^p} \leq  2p\left\|[Q]_{a, b}^{1/2}\right\|_{L^p}.
\end{equation}
Thus, if $[Q]_{a, b} \leq v$ almost surely for some deterministic $v > 0$, Markov's inequality gives \begin{equation}\label{eq:MarkovBDGBound}
    \Pp(Q_b^* > 2ep\sqrt{v}) \leq e^{-p}.
\end{equation}
For a complex continuous local martingale $Q$ with $Q_a = 0$ whose real and imaginary parts each have quadratic variation increment at most $v$ on $[a, b]$, Minkowski's inequality gives $\|Q_b^*\|_{L^p} \leq 4p\sqrt{v}$, and hence \eqref{eq:MarkovBDGBound} holds with $2e$ replaced by $4e$. 

We now establish the resolvent bound \eqref{eq:RigidityCoarseBoundDeltaandQ} which will be used to estimate the drift integral in \eqref{eq:LogFieldTerminalTime} over $[0,t_1]$. It also bounds $\delta_t(z_t)$ and $\partial_z\delta_t(z_t)$ at $t=t_1$, providing the initial estimates for our analysis of their evolution on $[t_1,T]$.
\begin{proposition}[Uniform rigidity]\label{prop:rigidityEvent} For a sufficiently large fixed choice of $A$ in \eqref{eq:intermediateHeightParams}, there exist constants $c, C > 0$ and an event $\Omega_N$ with 
\begin{equation}\label{eq:RigidityFailureProb}
    \Pp(\Omega_N^c) \leq e^{-r_N^c},
\end{equation}on which the estimates \begin{equation}\label{eq:EigenvalueRigidity}
     |\lambda_i(t) - \gamma_i| \leq Cr_NN^{-2/3}\min(i, N+1-i)^{-1/3} \quad (1 \leq i \leq N), 
\end{equation}and \begin{equation}\label{eq:countingRigidity} 
    \sup_{E \in \R}|\#\{i : \lambda_i(t) \leq E\} - NF_\scal(E)| \leq Cr_N 
\end{equation}hold simultaneously for all $t \in [0, T]$. On the same event, for every $z = E + \mathrm{i}\eta \in \C_+$ and $0 \leq t \leq T$, \begin{equation}\label{eq:RigidityCoarseBoundDeltaandQ}
    |\delta_t(z)| + \eta|\partial_z\delta_t(z)|\leq \frac{Cr_N}{N\eta}. 
\end{equation}\end{proposition}
\begin{proof} 

We first apply \cite[Theorem 2.2]{EYY} to establish \eqref{eq:EigenvalueRigidity} and \eqref{eq:countingRigidity} at each fixed time $t \in [0,T]$. We then extend these estimates to hold simultaneously for all $t \in [0,T]$ by controlling the matrix increments between points of a deterministic time mesh. Finally, we deduce \eqref{eq:RigidityCoarseBoundDeltaandQ} from \eqref{eq:countingRigidity} by Stieltjes integration by parts. \vspace{2mm}

    {\it Fixed-time rigidity. }At each deterministic time, \eqref{eq:DynamicalModel} gives independent centered upper-triangular entries with off-diagonal variance $1/N$ and diagonal variance $(1 + e^{-t})/N$. To apply \cite[Theorem 2.2]{EYY}, set $$w_t := 1 + \frac{e^{-t}}{N}, \qquad \widehat{H}_t := w_t^{-1/2}H_t.$$
    The variance matrix of $\widehat{H}_t$ is $$(\E|\hat{h}_{ij}|^2)_{i, j} = \frac{1}{Nw_t}\mathbf{1}\mathbf{1}^\top + \frac{e^{-t}}{Nw_t}\operatorname{Id},$$ where $\mathbf{1}$ is the vector of all ones. Its row sums equal one, its eigenvalue one is simple, and all other eigenvalues are bounded above by $1/N$. Together with the uniform Gaussian tail bounds, these properties verify the hypotheses of  \cite[Theorem 2.2]{EYY} with constants independent of $t.$

    We take $L = A\log\log N$, which lies in the admissible range of \cite[Theorem 2.2]{EYY} for all sufficiently large $N$, provided that $A$ is fixed and sufficiently large. \cite[Theorem 2.2]{EYY} then implies that \eqref{eq:EigenvalueRigidity} and \eqref{eq:countingRigidity} hold for $\widehat{H}_t$, for each fixed $t\geq0$ with 
    %The theorem then gives the rigidity and counting estimates with the factor  $(\log N)^L = r_N$, and 
    failure probability at most $C\exp(-cr_N^\phi)$ for some fixed $\phi > 0$. Since $r_N \to \infty$, this failure probability is at most $e^{-r_N^{c'}}$ for any fixed $0 < c' < \phi$ and all sufficiently large $N$.

    Returning from $\widehat{H}_t$ to $H_t$ changes the eigenvalue locations by $O(N^{-1})$ on the rigidity event. Since $F_\scal$ is globally Lipschitz and $w_t = 1 + O(N^{-1})$, we have $$N\left|F_\scal(E/\sqrt{w_t}) - F_\scal(E)\right| = O(1)$$uniformly for $E \in [-3, 3]$. These errors are absorbed by the constants in \eqref{eq:EigenvalueRigidity}-\eqref{eq:countingRigidity}. On the same event, the spectrum lies in $[-3, 3]$ for sufficiently large $N$, so the counting estimates extend to every $E \in \R$. \vspace{2mm}

    {\it Passage to all times.} Take a deterministic mesh of $[0, T]$ with spacing at most $\Delta:= N^{-12}$, including both endpoints and containing at most $CN^{12}$ points. Since $r_N^c/\log N \to \infty$, a union bound over the mesh points and the failure events for rigidity and the counting estimates shows that both estimates hold simultaneously at every mesh point, with failure probability at most $e^{-r_N^c}$ after reducing $c$. 

    Multiplying $H_t$ by $e^{t/2}$ cancels the drift term $-\frac{1}{2}H_t\dd t$ in \eqref{eq:OrnsteinUhlenbeck}. The resulting process $$Y_t:= e^{t/2}H_t = G_N + N^{-1/2}\int_0^te^{u/2}\dd B_u $$ is a (matrix-valued) martingale with respect to $(\mathcal{F}_t)$, so we can apply \eqref{eq:LpBDGBound} to its coordinates. For each real coordinate (that is, a diagonal entry or the real or imaginary part of an off-diagonal entry) of $Y_t - Y_s$, the quadratic variation accumulated over the mesh interval $[s, (s + \Delta) \wedge T]$ is at most $CN^{-13}$. Applying \eqref{eq:LpBDGBound} with $p = N$ to the real and imaginary parts, and then using the triangle inequality, gives $$\left\|\sup_{s \leq t \leq (s + \Delta)\wedge T}|(Y_t - Y_s)_{ij}|\right\|_{L^N} \leq CN\sqrt{N^{-13} } = CN^{-11/2}.$$
    By Markov's inequality, the probability that any fixed such increment exceeds $N^{-5}$ is at most $(CN^{-1/2})^N$. Taking a union bound over these events for all $O(N^2)$ real matrix coordinates and $O(N^{12})$ mesh intervals, we obtain a failure probability of at most $e^{-cN\log N}$. 

    Since the operator norm is bounded by the Frobenius norm, the preceding entrywise estimate implies $$\sup_{s \leq t \leq (s + \Delta)\wedge T}\|Y_t - Y_s\|_{\text{op}} \leq CN^{-4}$$simultaneously on every mesh interval. At a mesh point, rigidity gives $\|H_s\|_{\text{op}} \leq 3$. The identity$$H_t - H_s = (e^{-(t -s)/2 } - 1)H_s + e^{-t/2}(Y_t - Y_s)$$therefore yields \begin{equation}\label{eq:OperatorNormBoundOnRigidity}
        \sup_{s\leq t \leq (s + \Delta)\wedge T}\|H_t - H_s\|_{\text{op}} \leq CN^{-4}. 
    \end{equation}
    Weyl's eigenvalue perturbation inequality now gives $$\max_{1 \leq i \leq N}|\lambda_i(t) - \lambda_i(s)| \leq \|H_t - H_s\|_{\text{op}} \leq CN^{-4}. $$Together with rigidity at $s$, this proves \eqref{eq:EigenvalueRigidity} at every time.  With $\varepsilon = CN^{-4}$, it also gives $$\mathcal{N}_N(E - \varepsilon ; s) \leq \mathcal{N}_N(E ;  t) \leq \mathcal{N}_N(E + \varepsilon ; s).$$
    Since $F_\scal$ is globally Lipschitz, \eqref{eq:countingRigidity} therefore extends to every time, with an additional error of at most $CN\varepsilon = O(N^{-3})$. Define $\Omega_N$ as the intersection of the event on which \eqref{eq:EigenvalueRigidity}-\eqref{eq:countingRigidity} hold at every mesh point, and the event on which all the preceding entrywise increment bounds hold. The failure probabilities prove \eqref{eq:RigidityFailureProb}.\vspace{2mm}

    {\it Resolvent bounds. } On $\Omega_N$, the difference between the empirical and semicircle distribution functions, $$F_t(x) := N^{-1}\mathcal{N}_N(x ; t) - F_\scal(x),$$ satisfies $\|F_t\|_\infty \leq Cr_N/N$. Stieltjes integration by parts gives $$\delta_t(z) = \int_\R\frac{F_t(x)}{(x - z)^2}\dd x, \qquad \partial_z\delta_t(z) = 2 \int_\R\frac{F_t(x)}{(x - z)^3}\dd x.$$
    The boundary term vanishes because both measures have total mass one. Since $$\int_\R|x-z|^{-2}\dd x = \frac{\pi}{\eta}, \qquad \int_\R|x-z|^{-3}\dd x = \frac{2}{\eta^2},$$we obtain \eqref{eq:RigidityCoarseBoundDeltaandQ}. 
\end{proof}
\begin{corollary}
\label{cor:IncomingError} On $\Omega_N$, \begin{equation}\label{eq:IncomingError}
    |\delta_{t_1}| \leq \frac{Cr_N}{Nh_N}, \qquad |q_{t_1}| \leq \frac{Cr_N}{Nh_N^2}, \qquad N\int_0^{t_1}|\delta_t|^2\dd t \leq Cr_N^{-2}. 
\end{equation}
\end{corollary}
\begin{proof}
    The first two estimates follow from \eqref{eq:RigidityCoarseBoundDeltaandQ} at height $\eta = h_N$. For the third, the characteristic integral inequalities from Lemma~\ref{lem:CharacteristicGeometry} imply $$N\int_0^{t_1}|\delta_t|^2\dd t \leq \frac{Cr_N^2}{N}\int_0^{t_1}y_t^{-2}\dd t \leq \frac{Cr_N^2}{Nh_N} = Cr_N^{-2}. $$
\end{proof}

\section{Fine control on the terminal interval}\label{sec:TerminalInt} 
 The coarse estimate \eqref{eq:RigidityCoarseBoundDeltaandQ} is sufficient on $[0, t_1]$, where $y_t \geq h_N$. On $[t_1, T]$, we use the bounds for $\delta_{t_1}$ and $q_{t_1}$ from Corollary~\ref{cor:IncomingError} as initial estimates and obtain sharper bounds from their evolution equations. We first control the martingale increments up to a stopping time, using the BDG inequality on intervals where the characteristic height decreases by a factor of two. We then choose the stopping time to be the first time that either $|\delta_t|$ or $|q_t|$ reaches a prescribed threshold, and use their evolution equations to show that the stopping does not occur before $T$. This follows the approach of Huang and Landon \cite[Section 3.2]{HuangLandon}.

 \begin{lemma}\label{lem:DyadicStoppedEstimate}
  Let $\sigma$ be a stopping time with $t_1 \leq \sigma \leq T$, and suppose that for a fixed $C_* > 0$, $$\operatorname{Im}m_t(z_t) \leq C_*$$for all $t_1 \leq t < \sigma$ almost surely. For every $p \geq 2$, there is an event $\mathcal{H}_N(p)$ satisfying \begin{equation}\label{eq:dyadicEventFailure}
    \Pp(\mathcal{H}_N(p)^c) \leq C(1 + \log(h_N/\eta_*))e^{-p}
\end{equation}on which, simultaneously for $t_1 \leq t \leq T$ and $k = 0, 1$, \begin{equation}\label{eq:martingaleDifference}
    \left|M_{t \wedge \sigma}^{(k)} - M_{t_1}^{(k)}\right| \leq \frac{Cp}{Ny_{t \wedge \sigma}^{k + 1}}.
\end{equation}The constants depend only on $C_*$ and the characteristic bounds, and are independent of $p, N$, and $\sigma$.
\end{lemma}
\begin{proof}
    Let $a_0 > a_1 > \cdots > a_J$ be a decreasing sequence of positive numbers satisfying $$a_0 = h_N, \qquad a_J = \eta_*, \qquad \frac{a_{j - 1}}{a_j} \leq 2 \quad (1 \leq j \leq J),$$and $J \leq C(1 + \log(h_N/\eta_*))$. Define the deterministic times $\tau_j$ by $y_{\tau_j} = a_j$, so that $\tau_0 = t_1$ and $\tau_J = T$.

    For either real component $Q^{(k)}$ of $M_{s\wedge \sigma}^{(k)} - M_{t_1}^{(k)}$, the estimates \eqref{eq:QuadVariationBoundDeltaMG}-\eqref{eq:QuadVariationBoundqMG} and the stopping assumption give $$\dd \langle Q^{(k)}\rangle_u \leq \frac{C\mathbf{1}_{\{u <\sigma\}}}{N^2y_u^{2k + 3}}\dd u. $$
    The characteristic integral inequalities in Lemma~\ref{lem:CharacteristicGeometry} therefore imply $$[Q^{(k)}]_{t_1, \tau_j} \leq \frac{C}{N^2}\int_{t_1}^{\tau_j}y_u^{-(2k + 3)}\dd u \leq \frac{C}{N^2a_j^{2k + 2}}.$$Applying the complex version of \eqref{eq:LpBDGBound} and then Markov's inequality, we obtain $$\Pp\left(\sup_{t_1 \leq s \leq \tau_j}|M_{s\wedge\sigma}^{(k)} - M_{t_1}^{(k)}| > \frac{Cp}{Na_j^{k+1}}\right)\leq e^{-p}, $$after absorbing $e$ into the constant. A union bound over $j$ and $k$ defines $\mathcal{H}_N(p)$ and proves \eqref{eq:dyadicEventFailure}. 

    On this event, fix $t$ and write $u = t \wedge \sigma$. If $u = t_1$, the increment is zero. Otherwise, choose $j$ with $\tau_{j - 1} \leq u \leq \tau_j$. Since $a_j \leq y_u \leq 2a_j$, $$|M_u^{(k)} - M_{t_1}^{(k)}|\leq \frac{Cp}{Na_j^{k + 1}} \leq \frac{2^{k + 1}Cp}{Ny_u^{k + 1}}.$$Absorbing $2^{k + 1} \leq 4$ into the constant proves \eqref{eq:martingaleDifference}. 
\end{proof}
For the remainder of the proof, we take \begin{equation}\label{eq:pDefinition}
    p_N := \lceil 16 \log r_N \rceil.
\end{equation}
This choice makes the failure probability in \eqref{eq:dyadicEventFailure} a negative power of $r_N$, while ensuring $p_N/(N\eta_*) \to 0$ as required by the estimates below. The specific comparisons used below are \begin{equation}\label{eq:pNScales}
    h_N \to 0, \qquad \frac{r_N}{Nh_N}(1 + \log(h_N/\eta_*)) \to 0, \quad \frac{p_N}{N\eta_*} \to 0, \quad \frac{p_N^2}{N\eta_*} \to 0.
\end{equation}
These follow from $r_N/(Nh_N) = r_N^{-3}, p_N = O((\log \log N)^2)$, and $N\eta_* = (\log N)^{1/4}$. 
\begin{proposition}[Resolvent bounds on the terminal interval]\label{prop:LocalLawShortScales} There is an event $\mathcal{G}_N \subseteq \Omega_N$ satisfying \begin{equation}\label{eq:LocalLawFailureProb}
    \Pp(\Omega_N \setminus \mathcal{G}_N)\leq r_N^{-12},
\end{equation}on which, uniformly for $t_1 \leq t \leq T$, \begin{equation}\label{eq:localLawShortScale}
    |\delta_t| \leq C\left(\frac{r_N}{Nh_N} + \frac{\log r_N}{Ny_t}\right), \qquad |q_t|\leq C\left(\frac{r_N}{Nh_Ny_t} + \frac{\log r_N}{Ny_t^2}\right).
\end{equation}
\end{proposition}
\begin{proof}
    Write $p = p_N$, and define the thresholds \begin{equation}\label{eq:thresholds}
    U_t = K\left(\frac{r_N}{Nh_N} + \frac{p}{Ny_t}\right), \quad V_t = K^2\left(\frac{r_N}{Nh_Ny_t} + \frac{p}{Ny_t^2}\right),\end{equation} where $K$ is a fixed constant to be chosen below. The additional factor $K$ in $V_t$ makes the contribution of $\delta_t\partial_z^2m_t(z_t)$ in \eqref{eq:qSDE} small relative to the proposed bound for $q_t$. Define the stopping time $$\sigma = T \wedge \inf\{t \in [t_1, T] : |\delta_t| \geq U_t \text{ or }|q_t| \geq V_t\}.$$
    For each fixed $K$, \eqref{eq:pNScales} implies $U_T \leq 1$ for large $N$. Since $U_t$ is nondecreasing, before $\sigma$ we have $$\operatorname{Im}m_t(z_t) \leq \operatorname{Im}m_\scal(z_t) + U_T \leq C_*,$$where $C_*$ is independent of $K$, provided $N$ is sufficiently large. We may therefore apply Lemma~\ref{lem:DyadicStoppedEstimate}. 

    We work on $\Omega_N \cap \mathcal{H}_N(p)$. By \eqref{eq:IncomingError}, the stopping time $\sigma$ is strictly greater than $t_1$ for sufficiently large fixed $K$. For $t_1 \leq t \leq \sigma$, integrating \eqref{eq:deltaSDE} and \eqref{eq:qSDE}, and applying the initial estimates in \eqref{eq:IncomingError} and the martingale bounds in \eqref{eq:martingaleDifference} gives \begin{align*}
    |\delta_t| &\leq C\left(\frac{r_N}{Nh_N} + \frac{p}{Ny_t}\right) + C\int_{t_1}^tU_u\dd u + \int_{t_1}^tU_uV_u \dd u, \\
    |q_t| &\leq Cy_t^{-1}\left(\frac{r_N}{Nh_N} + \frac{p}{Ny_t}\right) + C\int_{t_1}^tV_u\dd u + \int_{t_1}^tV_u^2\dd u + C\int_{t_1}^t\frac{U_u}{y_u^2}\dd u.
\end{align*}Here, we used $|q_{t_1}| \leq Cr_N/(Nh_Ny_t)$, which follows from \eqref{eq:IncomingError} and $y_t \leq h_N$, and the bound  $|\partial_z^2m_u(z_u)| \leq 2C_*/y_u^2$ from \eqref{eq:DeterministicResolventInequalities} before $\sigma$. Using $(a + b)^2 \leq 2(a^2 + b^2)$ and
\eqref{eq:characteristicIntIneq}, we obtain
\begin{align*}
    \int_{t_1}^t U_uV_u\dd u
    &\leq CK^3\left[
        \left(\frac{r_N}{Nh_N}\right)^2
        \log\frac{h_N}{y_t}
        +\frac{p^2}{N^2y_t^2}
    \right],\\
    \int_{t_1}^t V_u^2\dd u
    &\leq CK^4\left[
        \left(\frac{r_N}{Nh_N}\right)^2
        \frac{1}{y_t}
        +\frac{p^2}{N^2y_t^3}
    \right],\\
    \int_{t_1}^t\frac{U_u}{y_u^2}\dd u
    &\leq CK\left[
        \frac{r_N}{Nh_Ny_t}
        +\frac{p}{Ny_t^2}
    \right].
\end{align*}
The integrals $\int_{t_1}^tU_u\dd u$ and $\int_{t_1}^tV_u\dd u$ are bounded by $Ch_NU_t$ and $Ch_NV_t$, respectively, because $T - t_1 \leq Ch_N$ and the thresholds are nondecreasing. Dividing by the thresholds and using $y_t \geq \eta_*$ yields \begin{align*}
    \frac{|\delta_t|}{U_t} &\leq \frac{C}{K} + Ch_N + CK^2\left[\frac{r_N}{Nh_N}\log\frac{h_N}{\eta_*} + \frac{p}{N\eta_*}\right], \\
    \frac{|q_t|}{V_t} &\leq \frac{C}{K^2} + \frac{C}{K} + Ch_N + CK^2\left[\frac{r_N}{Nh_N} + \frac{p}{N\eta_*}\right].
\end{align*}

Here, the constants $C$ are independent of $K$. Choose $K$ sufficiently large, and then $N$ sufficiently large, using \eqref{eq:pNScales}, so that both right-hand sides are at most $1/2$. By continuity, we conclude that $\sigma = T$ on this event, and the improved bounds hold at $T$. Finally, $$\Pp(\mathcal{H}_N(p_N)^c) \leq C(1 + \log r_N)r_N^{-16} \leq r_N^{-12}$$for large $N$. Setting $\mathcal{G}_N := \Omega_N \cap \mathcal{H}_N(p_N)$ proves \eqref{eq:LocalLawFailureProb}. Since $K$ is fixed and $p_N = O(\log r_N)$, the estimates obtained above imply \eqref{eq:localLawShortScale}. 
\end{proof}
The event $\Omega_N$ has overwhelming probability, whereas the bound for $\mathcal{G}_N$ only gives failure probability of order at most $r_N^{-12}$. This weaker bound is sufficient for the $L^2$ estimates in Propositions~\ref{prop:secondMomentSup} and \ref{prop:L2DriftSmall}, since the coarse estimate \eqref{eq:RigidityCoarseBoundDeltaandQ} remains available on $\Omega_N \setminus \mathcal{G}_N$. 
\section{Second moments and the drift estimate}\label{sec:SecondMomentBounds}
Our goal in this section is to show that the quadratic drift term $$-\frac{N}{2}\int_0^T\delta_t^2\dd t$$in \eqref{eq:LogFieldTerminalTime} tends to zero in $L^2$. We first prove a uniform second-moment estimate for the resolvent along the characteristic, which will also be used to control the martingale term. 
\begin{proposition}\label{prop:secondMomentSup}We have \begin{equation}\label{eq:SupL2NormDelta}
    \left\|\sup_{0 \leq t \leq T}|m_t(z_t) - m_\scal(z_t)|\right\|_{L^2} \leq C\left(\frac{r_N}{Nh_N} + \frac{\log r_N}{N\eta_*}\right) = o(1).
\end{equation}
Consequently, \begin{equation}\label{eq:L1SupImaginaryPart}
    \E \sup_{0 \leq t \leq T}|m_t(z_t)|^2 \leq C, \qquad \E\sup_{0 \leq t \leq T}\operatorname{Im}m_t(z_t) \leq C.
\end{equation}
\end{proposition}
\begin{proof}
    Let $S_N = \sup_{t \leq T}|\delta_t|$. On $\mathcal{G}_N$, the coarse estimate \eqref{eq:RigidityCoarseBoundDeltaandQ} before $t_1$ and \eqref{eq:localLawShortScale} afterward give $$S_N \leq C\left(\frac{r_N}{Nh_N} + \frac{\log r_N}{N\eta_*}\right).$$
    On $\Omega_N \setminus \mathcal{G}_N$, \eqref{eq:RigidityCoarseBoundDeltaandQ} gives $S_N \leq Cr_N/(N\eta_*)$. On $\Omega_N^c$, the deterministic bound $|m_t(z_t)| \leq y_t^{-1}$ and boundedness of $m_\scal$ give $S_N \leq C/\eta_*$. Splitting over these three disjoint events, and using \eqref{eq:RigidityFailureProb} and \eqref{eq:LocalLawFailureProb}, we obtain $$\|S_N\|_{L^2} \leq C\left(\frac{r_N}{Nh_N} + \frac{\log r_N}{N\eta_*} \right) + \frac{Cr_N^{-5}}{N\eta_*} + \frac{C}{\eta_*}e^{-r_N^c/2}.$$ 
    Since $r_N^{-5} = o(\log r_N)$, and $$\frac{1}{\eta_*}e^{-r_N^c/2} = o\left(\frac{r_N}{Nh_N}\right),$$the last two terms are negligible compared with the first bound. This proves \eqref{eq:SupL2NormDelta}. Boundedness of $m_\scal(z_t)$ and Cauchy-Schwarz prove \eqref{eq:L1SupImaginaryPart}. 
\end{proof}
\begin{proposition}[Drift estimate]\label{prop:L2DriftSmall}
Let $$D_N := N \int_0^T|\delta_t|^2\dd t.$$
Then \begin{equation}\label{eq:DriftTermSmall}
    \|D_N\|_{L^2} \leq C\left(r_N^{-2} + \frac{(\log r_N)^2}{N\eta_*}\right) = o(1).
\end{equation}In particular, $\|\tfrac{N}{2}\int_0^T\delta_t^2\dd t\|_{L^2} = o(1).$
\end{proposition}
\begin{proof}
    On $\mathcal{G}_N$, Corollary~\ref{cor:IncomingError} controls the contribution from $[0, t_{1}]$. Combining this with \eqref{eq:localLawShortScale}, we obtain\begin{align*}
        D_N &\leq Cr_N^{-2} + CN\left(\frac{r_N}{Nh_N}\right)^2h_N + \frac{C(\log r_N)^2}{N}\int_{t_1}^Ty_t^{-2}\dd t \\
        &\leq C\left(r_N^{-2} + \frac{(\log r_N)^2}{N\eta_*}\right). 
    \end{align*}
    On $\Omega_N$, the coarse estimate \eqref{eq:RigidityCoarseBoundDeltaandQ} gives $$D_N \leq \frac{Cr_N^2}{N}\int_0^Ty_t^{-2}\dd t \leq \frac{Cr_N^2}{N\eta_*}. $$
    On $\Omega^c_N$, the bound $|\delta_t|\leq C/y_t$ implies $D_N \leq CN/\eta_*$. Thus, the same decomposition as in the previous proposition yields $$\|D_N\|_{L^2} \leq C\left(r_N^{-2} + \frac{(\log r_N)^2}{N\eta_*}\right) + \frac{Cr_N^{-4}}{N\eta_*} + \frac{CN}{\eta_*}e^{-r_N^c/2}.$$
    Again, the last two terms are negligible and can be absorbed into the first. The result follows since $$\frac{(\log r_N)^2}{N\eta_*} = \frac{A^2(\log \log N)^4}{(\log N)^{1/4}} \to 0. $$
\end{proof}
\section{Gaussian approximation and evaluation of its variance}\label{sec:GaussianSection}
We now approximate the martingale in \eqref{eq:LogFieldTerminalTime} by a Gaussian random variable and compute its variance. Split the integral in \eqref{eq:LogFieldTerminalTime} at $t_1$ by writing $$M_N^{\text{main}} := \frac{1}{\sqrt{N}}\sum_i\int_0^{t_1}\frac{\dd B_i(t)}{\lambda_i(t) - z_t}, \qquad M_N^{\text{tail}} := \frac{1}{\sqrt{N}}\sum_i\int_{t_1}^{T}\frac{\dd B_i(t)}{\lambda_i(t) - z_t},$$so that $M_N =  M_N^{\text{main}} + M_N^{\text{tail}}$. 
On $[0, t_1]$, replace the eigenvalues by their semicircle quantiles, and define \begin{equation}\label{eq:ExactGaussian}
    \widetilde{M}_N := \frac{1}{\sqrt{N}}\sum_i\int_0^{t_1}\frac{\dd B_i(t)}{\gamma_i - z_t}, \qquad Z_N := \operatorname{Im}\widetilde{M}_N. 
\end{equation}Note that all integrands in \eqref{eq:ExactGaussian} are deterministic. By Lemma~\ref{lem:EigenvalueDBMandSDEs}, $Z_N$ is therefore a centered real Gaussian independent of $H_0$. Landon and Sosoe \cite[Section~7]{LandonSosoe} use the
same approximation of the stochastic integral by replacing
the eigenvalues with their semicircle quantiles, and
compute the variance of the resulting Gaussian variable. The following estimate will be used to bound
$M_N^{\mathrm{main}}-\widetilde M_N$ and to compute
$\operatorname{Var}(Z_N)$. 
\begin{lemma}\label{lem:QuantileIntegration}
    Let $f : \R \to \C$ be absolutely continuous, vanish at infinity, and satisfy $f'\in L^1(\R)$. Then \begin{equation}\label{eq:GammavsSemicircleIntegral}
        \left|\frac{1}{N}\sum_if(\gamma_i) - \int_\R f(x)\rho_\scal(x)\dd x\right| \leq \frac{1}{N}\int_\R|f'(x)|\dd x. 
    \end{equation}In particular, whenever $N\eta \geq 1$, \begin{align}\label{eq:FourthMomentGammaIntegralBound}
        \frac{1}{N}\sum_i|\gamma_i - E - \mathrm{i}\eta|^{-4} \leq C\eta^{-3}.
    \end{align}
\end{lemma}
\begin{proof}
    The distribution functions of $N^{-1}\sum_i\delta_{\gamma_i}$ and the semicircle law differ by at most $N^{-1}$, so integration by parts proves \eqref{eq:GammavsSemicircleIntegral}. For $f(x) = ((x - E)^2 + \eta^2)^{-2}$, we have $$\int_\R|f'(x)|\dd x = 2\eta^{-4}, \qquad \int_\R f(x)\rho_\scal(x)\dd x \leq C\eta^{-3}.$$Since $2/(N\eta^4)$ is absorbed into $C\eta^{-3}$ when $N\eta \geq 1$, this proves \eqref{eq:FourthMomentGammaIntegralBound}.
\end{proof}
\begin{proposition}[Gaussian approximation and  variance]\label{prop:MartingaleApproximation}The stochastic integrals defined above satisfy \begin{align}
    \|M_N^{\text{\normalfont main}} - \widetilde{M}_N\|_{L^2} &= o(1), \label{eq:MainMartingaleL2} \\
    \E|M_N^{\text{\normalfont tail}}|^2 &\leq C\log(h_N/\eta_*) = O(\log r_N). \label{eq:L2TailMartingaleSmall}
\end{align}
The variance of $Z_N$ is \begin{equation}\label{eq:VarianceExact}
    v_N := \operatorname{Var}(Z_N) = \frac{1}{N}\sum_{i = 1}^N\int_0^{t_1}\frac{y_t^2}{((\gamma_i - x_t)^2 + y_t^2)^2}\dd t. 
\end{equation}
It has the asymptotic expansion \begin{equation}\label{eq:varianceAsymptoticComputed}
    v_N = \frac{1}{2}\log\frac{y_0}{h_N} + O\left(T + \frac{1}{Nh_N}\right),
\end{equation}and hence \begin{align}
    v_N &= \frac{\vartheta}{2}\log N - 2\log r_N + O(1) = \frac{\vartheta}{2}\log N + o(\log N). \label{eq:varianceAsymptoticTheta}
\end{align} 
\end{proposition}
\begin{proof} We first estimate $M_N^{\text{main}}-\widetilde M_N$ using eigenvalue rigidity on $\Omega_N$ and a deterministic
bound on $\Omega_N^c$. We then bound $M_N^{\text{tail}}$ using Proposition~\ref{prop:secondMomentSup}. Finally, we compute the variance of $Z_N$ by approximating the sum over semicircle quantiles by an integral against $\rho_\scal$ and integrating along the characteristic.\vspace{2mm}

    {\it The main martingale.} Fix a compact interval $J \subset\subset (-2, 2)$ such that all $x_t$ lie in its interior at a fixed positive distance from its boundary. On $\Omega_N$, the rigidity estimate \eqref{eq:EigenvalueRigidity} gives $|\lambda_i(t) - \gamma_i| \leq Cr_N/N$ whenever $\gamma_i \in J$, since $J$ is a compact subset of $(-2, 2)$. For $t \leq t_1$, the bound $y_t \geq h_N$ and \eqref{eq:RigidityVsHeightRatio} imply $$|\lambda_i(t) - z_t| \geq \tfrac{1}{2}|\gamma_i - z_t|$$for large $N$. Thus \begin{equation}\label{eq:BulkRigidityReplacement}
        \left|\frac{1}{\lambda_i(t) - z_t} - \frac{1}{\gamma_i - z_t}\right|^2 = \frac{|\lambda_i(t) - \gamma_i|^2}{|\lambda_i(t) - z_t|^2|\gamma_i - z_t|^2} \leq \frac{Cr_N^2}{N^2}|\gamma_i - z_t|^{-4}. 
    \end{equation}
    For $\gamma_i \not\in J$, both denominators are bounded below by a fixed constant, since $x_t$ is separated from $J^c$ and the largest rigidity error $Cr_NN^{-2/3}$ tends to zero. Moreover, $$\frac{1}{N}\sum_i|\lambda_i(t) - \gamma_i|^2 \leq Cr_N^2N^{-7/3}\sum_i\min(i, N + 1 - i)^{-2/3} \leq \frac{Cr_N^2}{N^2}.$$ Combining these bounds with \eqref{eq:BulkRigidityReplacement} and Lemma~\ref{lem:QuantileIntegration}, we obtain, on $\Omega_N$,  \begin{equation}\label{eq:replacementError}
        \frac{1}{N}\sum_i\left|\frac{1}{\lambda_i(t) - z_t} - \frac{1}{\gamma_i - z_t}\right|^2 \leq \frac{Cr_N^2}{N^2}(y_t^{-3} + 1). 
    \end{equation}
    For every configuration (without rigidity), the left-hand side is at most $4/y_t^2$.

    Apply Itô's isometry, and split the resulting expectation over $\Omega_N$ and $\Omega_N^c$. By \eqref{eq:characteristicIntIneq} with powers three and two, \begin{align*}
        \E|M_N^{\text{main}} - \widetilde{M}_N|^2 &= \int_0^{t_1}\E\left[\frac{1}{N}\sum_i \left|\frac{1}{\lambda_i(t) - z_t} - \frac{1}{\gamma_i - z_t}\right|^2\right]\dd t \\
        &\leq C\left[\left(\frac{r_N}{Nh_N}\right)^2 + \frac{Tr_N^2}{N^2} + \frac{e^{-r_N^c}}{h_N}\right] = o(1).
    \end{align*}The three terms vanish because $r_N/(Nh_N) = r_N^{-3}$, $r_N = O(N^\varepsilon)$ for every fixed $\varepsilon > 0$, and $r_N^c \gg \log N$. 

    {\it The terminal martingale. } Itô's isometry and \eqref{eq:DeterministicResolventInequalities} give$$\E|M_N^{\text{tail}}|^2 = \int_{t_1}^T\E\left[\frac{1}{N}\sum_i|\lambda_i(t) - z_t|^{-2}\right]\dd t = \int_{t_1}^T\frac{\E \operatorname{Im}m_t(z_t)}{y_t}\dd t. $$
    Proposition~\ref{prop:secondMomentSup} bounds the expected imaginary part uniformly, so \eqref{eq:characteristicIntIneq} proves \eqref{eq:L2TailMartingaleSmall}. \vspace{2mm}

    {\it The variance. } Taking imaginary parts in \eqref{eq:ExactGaussian} and applying Itô's isometry proves \eqref{eq:VarianceExact}. For $$f_{E, \eta}(x) := \frac{\eta^2}{((x - E)^2 + \eta^2)^2},$$we have $\int_\R |f'_{E, \eta}(x)|\dd x = 2/\eta^2$. Lemma~\ref{lem:QuantileIntegration} and \eqref{eq:characteristicIntIneq} therefore imply \begin{equation}\label{eq:varianceDoubleIntegral}
        v_N = \int_0^{t_1}\int_\R \frac{y_t^2\rho_\scal(x)}{((x - x_t)^2 + y_t^2)^2}\dd x\dd t + O\left(\frac{1}{Nh_N}\right). 
    \end{equation}
    The identity $$\frac{\eta^2}{((x - E)^2 + \eta^2)^2} = \frac{1}{2\eta}\operatorname{Im}\frac{1}{x - E - \mathrm{i}\eta} - \frac{1}{2}\operatorname{Re}\frac{1}{(x - E - \mathrm{i}\eta)^2}$$ expresses the inner integral as $$\frac{\operatorname{Im}m_\scal(E + \mathrm{i}\eta)}{2\eta} - \frac{1}{2}\operatorname{Re}m_\scal'(E + \mathrm{i}\eta). $$
    Using $\operatorname{Im}b(z_t) = \operatorname{Im}m_\scal(z_t) + y_t/2$, we obtain \begin{align*}
        v_N &= \frac{1}{2}\int_0^{t_1}\frac{\operatorname{Im}b(z_t)}{y_t}\dd t - \frac{t_1}{4} - \frac{1}{2}\int_0^{t_1}\operatorname{Re}m_\scal'(z_t)\dd t + O\left(\frac{1}{Nh_N}\right) \\
        &= \frac{1}{2}\log\frac{y_0}{h_N} + O\left(T + \frac{1}{Nh_N}\right).
    \end{align*}The last equality follows from \eqref{eq:characteristicIntEq} and boundedness of $m_\scal'$ along the characteristic. Since $Nh_N = r_N^4 \to \infty$, the error is $o(1)$. Finally, $$\log y_0 = -(1 - \vartheta)\log N + O(1), \qquad \log h_N = 4\log r_N - \log N,$$which proves \eqref{eq:varianceAsymptoticTheta}. 
\end{proof}
\begin{corollary}[Log-field decomposition]\label{cor:LogFieldDecomp}
Let $X_N := \operatorname{Im}L_0(z_0)$. There is a remainder $R_N$ such that \begin{equation}\label{eq:LogFieldDecomposition}
    \operatorname{Im}L_T(z_T) = X_N + Z_N + R_N, \qquad \|R_N\|_{L^2} = O(\sqrt{\log r_N}) + o(1) = o(\sqrt{\log N}).
\end{equation}The variables $X_N$ and $Z_N$ are independent, and $Z_N \sim \operatorname{Normal}(0, v_N)$ with $v_N$ as evaluated in Proposition~\ref{prop:MartingaleApproximation}. 
\end{corollary}
\begin{proof}
    Take imaginary parts in \eqref{eq:LogFieldTerminalTime}, split the martingale at $t_1$, and replace its main part by $\widetilde{M}_N$. Proposition~\ref{prop:L2DriftSmall} and \eqref{eq:MainMartingaleL2} give $o(1)$ errors in $L^2$, while \eqref{eq:L2TailMartingaleSmall} contributes $O(\sqrt{\log r_N})$. The integrands in \eqref{eq:ExactGaussian} are deterministic and the Brownian motions are independent of $H_0$, so $Z_N$ is independent of $X_N = \operatorname{Im}L_0(z_0).$
\end{proof}
\section{The initial GOE logarithmic field}\label{sec:InitialGOEField}
We next determine the fluctuations of $\operatorname{Im}L_0(z_0)$. We record the following specialization of \cite[Theorem B.1]{BLZ} to the GOE logarithmic field at the mesoscopic heights needed here, and apply the resulting estimate at the initial characteristic point $z_0$. 
\begin{lemma}\label{lem:MesoscopicGOECLT}
Let $z_N = E_N + \mathrm{i}\eta_N$ be deterministic, with $E_N$ in a fixed compact subinterval of $(-2, 2)$, $\eta_N \to 0$, and $\eta_N \geq (\log N)^{1000}/N$. For the GOE normalization in Section~\ref{sec:Intro}, define $X_N(z_N) := \operatorname{Im}L_0(z_N)$. Then \begin{equation}\label{eq:GOECLTInputMeanVariance}
    \E X_N(z_N) = O(1), \qquad \operatorname{Var}(X_N(z_N)) = \log(1/\eta_N) + O(1),
\end{equation}and \begin{equation}\label{eq:GOECLTInputCLT}
    \frac{X_N(z_N)}{\sqrt{\log(1/\eta_N)}} \xrightarrow{\mathrm{d}} \operatorname{Normal}(0, 1).
\end{equation}More precisely, there are real deterministic coefficients $\mu_N, V_N^{\text{\normalfont GOE}}$ with \begin{equation}\label{eq:GOECLTInputCoefficients}
    \mu_N = O(1), \qquad V_N^{\text{\normalfont GOE}} = \log(1/\eta_N) + O(1),
\end{equation}such that, uniformly for $s \in \C$ with $|s|\leq 1/20$, \begin{equation}\label{eq:GOECLTInputCharacteristicFunction}
    \E e^{sX_N(z_N)} = \exp\left(s\mu_N + \frac{s^2}{2}V_N^{\text{\normalfont GOE}}\right)\left(1 + O((N\eta_N)^{-1/2})\right).
\end{equation}
\end{lemma}
\begin{proof}
    In the convention of \cite{BLZ}, our GOE is the $\beta = 1$ ensemble with potential $V(x) = x^2/2$, which satisfies the assumptions of  \cite[Theorem B.1]{BLZ}. Apply that theorem at $z_N$ with test function $$h(x) = s\operatorname{Im}\log(z_N  - x),$$corresponding to $p = 1, \beta = 1, \zeta_1 = 0$, and $\xi_1 = s$ in their notation. Note we take the lower height cutoff in \cite[Theorem B.1]{BLZ} to be $\eta_N$. 

    For $x \in \R$ and $z \in \C_+$, our logarithm conventions are related by $\operatorname{Log}_-(x - z) = \log(z - x) - \mathrm{i}\pi$. The constant cancels after centering by the deterministic semicircle term, so the theorem gives \eqref{eq:GOECLTInputCharacteristicFunction}. The linear coefficient $\mu_N$ satisfies $\mu_N = O(1)$ by \cite[Lemma B.2]{BLZ}. Specializing the quadratic form in \cite[Equation (B.6)]{BLZ} to the semicircle law gives $$ V_N^{\text{\normalfont GOE}} = \log|1 - m_\scal(z_N)^2| - \log (1 - |m_\scal(z_N)|^2).$$
    On the other hand, the semicircle equation implies $$1 - |m_\scal(z_N)|^2 = \frac{\eta_N|m_\scal(z_N)|^2}{\operatorname{Im}m_\scal(z_N)}.$$Since $E_N$ remains in a fixed compact bulk interval, the quantities $|m_\scal(z_N)|, \operatorname{Im}m_\scal(z_N)$, and $|1 - m_\scal(z_N)^2|$ are uniformly bounded and bounded away from zero for large $N$. Consequently, $$V_N^{\text{\normalfont GOE}} = \log(1/\eta_N) + O(1).$$

    For fixed $N$, the bound $|X_N(z_N)| \leq \pi N$ implies that its moment generating function is entire. By \eqref{eq:GOECLTInputCharacteristicFunction}, the function $$e^{-s\mu_N - \frac{1}{2}s^2V_N^{\text{\normalfont GOE}}}\E e^{sX_N(z_N)} - 1$$is analytic and uniformly $O((N\eta_N)^{-1/2})$ on a fixed complex disk. Cauchy's estimates therefore give the same bound for its first two derivatives at zero. Differentiating the logarithm of the moment-generating function yields $$\E X_N(z_N) = \mu_N + O((N\eta_N)^{-1/2}), \qquad \operatorname{Var}(X_N(z_N)) = V_N^{\text{\normalfont GOE}} + O((N\eta_N)^{-1/2}),$$which proves \eqref{eq:GOECLTInputMeanVariance}. Finally, for fixed $u \in \R$, substitute $s = \mathrm{i}u/\sqrt{\log(1/\eta_N)}$ into \eqref{eq:GOECLTInputCharacteristicFunction}. By \eqref{eq:GOECLTInputCoefficients}, the characteristic function tends to $e^{-u^2/2}$. Lévy's continuity theorem therefore gives \eqref{eq:GOECLTInputCLT}. 
\end{proof}
At the initial characteristic point $z_0 = z_0^{(N)}$, Lemma~\ref{lem:CharacteristicGeometry} gives $x_0 = E_0 + O(T)$, and $Ny_0 \asymp N^\vartheta \gg (\log N)^{1000}$. Thus the lemma applies to $X_N = \operatorname{Im}L_0(z_0)$ and we obtain
\begin{equation}\label{eq:CLTForXN}
    \E X_N = O(1), \qquad \operatorname{Var}(X_N) = (1 - \vartheta)\log N + O(1), \qquad \frac{X_N}{\sqrt{(1 - \vartheta)\log N}}\xrightarrow{\mathrm{d}} \operatorname{Normal}(0, 1). 
\end{equation}
\section{Removal of regularization and proof of the main theorem}\label{sec:RemoveRegularization} 
We now relate the logarithmic field to the one-sided counting function. For $\eta > 0$, define \begin{gather*}
    q_\eta(u) :=  \frac{1}{2} - \frac{1}{\pi}\arctan(u/\eta), \\
    \mathcal{N}_{N, \eta}(E) := \sum_{i = 1}^N q_\eta(\lambda_i(T) - E), \qquad D_{N, \eta}(E):=N\int_\R q_\eta(x - E)\rho_\scal(x)\dd x. 
\end{gather*}
The branch convention in \eqref{eq:LogFieldDefinition} gives the exact identity \begin{equation}\label{eq:SmoothingCountIdentity}
    \mathcal{N}_{N, \eta}(E) - D_{N, \eta}(E) = -\frac{1}{\pi}\operatorname{Im}L_T(E + \mathrm{i}\eta).
\end{equation}
Here $\mathcal{N}_N(E) = \#\{i : \lambda_i(T) \leq E\}$, whose distribution agrees with that of the count in Theorem~\ref{thm:MainResult}. The following estimate controls the error in replacing this count by its regularization. 
\begin{lemma}\label{lem:DesmoothingEstimate}
For $E\in \R$ and $\eta > 0$, write $z = E + \mathrm{i}\eta$.  Then the deterministic bound \begin{equation}
\begin{aligned}
&\left|[\mathcal{N}_N(E) - NF_\scal(E)]
    - [\mathcal{N}_{N,\eta}(E) - D_{N,\eta}(E)]\right| \\
&\qquad \leq CN\eta\left(
    |\operatorname{Re}\delta_T(z)|
    + \operatorname{Im}m_T(z)
    + \operatorname{Im}m_\scal(z)
    \right)
\end{aligned}
\label{eq:desmoothingEstimate}
\end{equation}holds. 
Consequently, \begin{equation}\label{eq:DesmoothingEstimateOurScale}
    \|[\mathcal{N}_N(E_0) - NF_\scal(E_0)] - [\mathcal{N}_{N, \eta_*}(E_0) - D_{N, \eta_*}(E_0)]\|_{L^2} \leq C(\log N)^{1/4}. 
\end{equation}
\end{lemma}
\begin{proof}
    Following \cite[Section 3.2]{BLZ}, we correct the difference between the indicator and its regularization by setting $$a_\eta(u) := \mathbf{1}_{\{u \leq 0\}} - q_\eta(u) + \frac{\eta u}{\pi(u^2 + \eta^2)}.$$
    We first show that \begin{equation}\label{eq:aFnBound}
        |a_\eta(u)| \leq C\frac{\eta^2}{u^2 + \eta^2}
    \end{equation}for all $u \in \R$. For $|u|\leq \eta$, this follows from the uniform boundedness of $a_\eta$, and the bound $\eta^2/(u^2 + \eta^2) \geq 1/2$. For $|u| > \eta$, we have the identity $$\mathbf{1}_{\{u \leq 0\}} - q_\eta(u) = -\pi^{-1}\arctan(\eta/u).$$ With $y = \eta/u$, the first order terms in the Taylor expansions of $\arctan y$ and $y/(1 + y^2)$ cancel. Indeed, $$\arctan y - \frac{y}{1 + y^2} = \int_0^y\frac{2v^2}{(1 + v^2)^2}\dd v = O(|y|^3), \qquad |y| \leq 1.$$
    Since $|y|^3 \leq 2y^2/(1 + y^2)$ in this range, \eqref{eq:aFnBound} follows.

    Now let $\nu = N^{-1}\sum_i\delta_{\lambda_i(T)} - \rho_\scal(x)\dd x$. Then $$\int_\R\left[\mathbf{1}_{\{x \leq E\}} - q_\eta(x - E)\right]\dd \nu(x) = \int_\R a_\eta(x - E)\dd \nu(x) - \frac{\eta}{\pi}\operatorname{Re}\delta_T(z).$$
     Using \eqref{eq:aFnBound} and bounding the integrals against the empirical measure and the semicircle measure separately, we obtain \begin{align*}
        \left|\int_\R a_\eta(x - E)\dd \nu(x)\right| &\leq \frac{C\eta^2}{N}\sum_i \frac{1}{(\lambda_i(T) - E)^2 + \eta^2} + C\eta^2\int_\R \frac{\rho_\scal(x)}{(x - E)^2 + \eta^2}\dd x \\
        &= C\eta(\operatorname{Im}m_T(z) + \operatorname{Im}m_\scal(z)). 
    \end{align*}
    Multiplying by $N$ proves \eqref{eq:desmoothingEstimate}. At $z = z_T$, Proposition~\ref{prop:secondMomentSup} gives $\|\delta_T(z_T)\|_{L^2} = o(1)$ and $\|\operatorname{Im}m_T(z_T)\|_{L^2} \leq C$, while $\operatorname{Im}m_\scal(z_T)$ is bounded. Taking $L^2$ norms in \eqref{eq:desmoothingEstimate}, and using $N\eta_* = (\log N)^{1/4}$, proves \eqref{eq:DesmoothingEstimateOurScale}.
\end{proof}

\begin{proof}[Proof of Theorem~\ref{thm:MainResult}]
    {\it The log-field CLT. } Set $$\sigma_N^2 := (1 - \vartheta)\log N + v_N.$$
    By \eqref{eq:varianceAsymptoticTheta}, \begin{equation}\label{eq:sigmaSquared}
        \sigma_N^2 = (1 - \vartheta/2)\log N - 2\log r_N + O(1) \sim (1 - \vartheta/2)\log N.
    \end{equation}Now, \eqref{eq:CLTForXN} gives $$\frac{X_N}{\sqrt{(1 - \vartheta)\log N}} \xrightarrow{\dd }\operatorname{Normal}(0, 1),$$while $Z_N/\sqrt{v_N}$ is exactly standard normal and independent of $X_N$. Moreover, $$\frac{(1 - \vartheta)\log N}{\sigma_N^2} \to \frac{1 - \vartheta}{1 - \vartheta/2}, \quad \frac{v_N}{\sigma_N^2} \to \frac{\vartheta/2}{1 - \vartheta/2}.$$
    Independence implies that the characteristic function of $(X_N + Z_N)/\sigma_N$ factors, and the preceding limits show that it converges to $e^{-u^2/2}$ for every $u \in \R$. Lévy's continuity theorem therefore gives $(X_N + Z_N)/\sigma_N \xrightarrow{\mathrm{d}}\operatorname{Normal}(0, 1)$. Independence and \eqref{eq:CLTForXN} also give \begin{equation}\label{eq:VarianceXplusZ}
        \operatorname{Var}(X_N + Z_N) = \sigma_N^2 + O(1).
    \end{equation}
    
{\it Passage to the counting function. } Combining Corollary~\ref{cor:LogFieldDecomp}, \eqref{eq:SmoothingCountIdentity}, and Lemma~\ref{lem:DesmoothingEstimate}, we obtain a real remainder $\widetilde{R}_N$ such that \begin{equation}\label{eq:CountingSmoothedWithRemainder}
        \mathcal{N}_N(E_0) - NF_\scal(E_0) = -\frac{1}{\pi}(X_N + Z_N) + \widetilde{R}_N, \qquad \|\widetilde{R}_N\|_{L^2} \leq C(\log N)^{1/4}.
    \end{equation}
    Indeed, the three error terms contribute $C(\log N)^{1/4} + O(\sqrt{\log r_N}) + o(1)$, and $\sqrt{\log r_N} = o((\log N)^{1/4})$. Since $\sigma_N \asymp \sqrt{\log N}$, Slutsky's theorem and \eqref{eq:sigmaSquared} prove \eqref{eq:MainResult}. \vspace{2mm}

    {\it The mean and variance.} By \eqref{eq:CountingSmoothedWithRemainder} and Cauchy-Schwarz, \begin{align*}
        \left|\operatorname{Var}(\mathcal{N}_N(E_0)) - \frac{1}{\pi^2}\operatorname{Var}(X_N + Z_N)\right|& \\
        \leq \frac{2}{\pi}\sqrt{\operatorname{Var}(X_N + Z_N)}\left\|\widetilde{R}_N\right\|_{L^2}& + \left\|\widetilde{R}_N\right\|_{L^2}^2 \leq C(\log N)^{3/4}. 
    \end{align*}
    Equations~\eqref{eq:sigmaSquared}-\eqref{eq:VarianceXplusZ} and the fact that $\log r_N = o((\log N)^{3/4})$ prove \eqref{eq:VarianceAsymptotic}. Finally, taking expectations in \eqref{eq:CountingSmoothedWithRemainder} gives $$\E\mathcal{N}_N(E_0) - NF_\scal(E_0) = O((\log N)^{1/4}) = o(\sqrt{\log N}),$$which proves \eqref{eq:MeanAsymptotic}. 
\end{proof}

\noindent\textbf{Joint fluctuations.} To conclude this paper, we explain the extension mentioned in Remark~\ref{rem:MultipointFluctuations}. For the two energies
$E_1$ and $E_2$ considered there, let $z_t^{(j)} = x_t^{(j)} + \mathrm{i}y_t^{(j)}$ be the
characteristic ending at $E_j + \mathrm{i}\eta_*$, for $j = 1,2$. The estimates for the drift term, the Gaussian approximation, and the approximation of the counting function apply separately along each characteristic, with constants uniform for the energies in a fixed compact subinterval of $(-2,2)$. Thus the decomposition into an initial term $X_N^{(j)}$, a Gaussian term $Z_N^{(j)}$, and an $L^2$ remainder of order $o(\sqrt{\log N})$ holds at both energies. The initial GOE terms $X_N^{(j)} = \operatorname{Im}L_0(z_0^{(j)})$ also require no additional estimates, since \cite[Theorem B.1]{BLZ} already gives their joint Gaussian asymptotics. The Gaussian terms $Z_N^{(j)}$ are stochastic integrals with deterministic integrands and are taken with respect to the same Brownian motions, which are independent of $G_N$. They therefore form a Gaussian vector independent of the initial GOE terms. Thus, it only remains to compute the covariances of the initial GOE terms and the stochastic integrals. 

By independence and the $L^2$ bounds on the error terms, $$\operatorname{Cov}\left(\mathcal{N}_N(E_1), \mathcal{N}_N(E_2)\right) = \frac{1}{\pi^2}\left[\operatorname{Cov}(X_N^{(1)}, X_N^{(2)}) + \operatorname{Cov}(Z_N^{(1)}, Z_N^{(2)})\right] + o(\log N).$$
The first covariance can be obtained from the covariance formula in \cite[Appendix B]{BLZ}. For the second, let $t_1^{(j)}$ satisfy $y_{t_1^{(j)}}^{(j)} = h_N$. Then Itô's isometry gives $$\operatorname{Cov}(Z_N^{(1)}, Z_N^{(2)}) = \frac{1}{N}\sum_i\int_0^{t_1^{(1)} \wedge t_1^{(2)}}\frac{y_t^{(1)}y_t^{(2)}}{|\gamma_i - z_t^{(1)}|^2|\gamma_i - z_t^{(2)}|^2}\dd t. $$
Lemma~\ref{lem:QuantileIntegration} allows us to replace
the sum by an integral against the semicircle density.
For $|E_1-E_2|=(1+o(1))N^{-a}$ with $0<a<1$,
evaluating this integral and the initial GOE covariance
gives
\begin{equation}\label{eq:TwoPointCountingCovariance}
    \operatorname{Cov}\left(\mathcal N_N(E_1),
    \mathcal N_N(E_2)\right) = \frac{a + \min(a, 1 -\vartheta)}{2\pi^2}\log N
    + o(\log N).
\end{equation}
We omit the details of these calculations. Together with the one-point variance asymptotics, this identifies the covariance matrix of the joint Gaussian limit. The same argument applies to any fixed number of energies satisfying analogous pairwise separation assumptions, with the covariance computed for each pair.

\section*{Acknowledgements}
The author thanks Benjamin Landon for suggesting this problem and for helpful discussions throughout the work. \\

\noindent \textbf{AI Disclosure.} The author developed the mathematical results and proofs and wrote the initial draft without using generative AI. ChatGPT (OpenAI; GPT-6 Astra Pro) was used only afterward to assist with proofreading, suggest improvements to the exposition, and check for possible mathematical errors. The author takes responsibility for the final manuscript.


\begin{thebibliography}{99}

\bibitem{AGZ}
G.~W.~Anderson, A.~Guionnet, and O.~Zeitouni,
\emph{An Introduction to Random Matrices},
Cambridge Studies in Advanced Mathematics, vol.~118,
Cambridge University Press, Cambridge, 2010.

\bibitem{BourgadeExtremeGaps}
P.~Bourgade,
Extreme gaps between eigenvalues of Wigner matrices,
\emph{J. Eur. Math. Soc.} \textbf{24} (2022), no.~8, 2823--2873.
\href{https://doi.org/10.4171/JEMS/1141}{doi:10.4171/JEMS/1141}.

\bibitem{BLZ}
P.~Bourgade, P.~Lopatto, and O.~Zeitouni,
Optimal rigidity and maximum of the characteristic polynomial of Wigner matrices,
\emph{Geom. Funct. Anal.} \textbf{35} (2025), no.~1, 161--253.
\href{https://doi.org/10.1007/s00039-025-00701-5}{doi:10.1007/s00039-025-00701-5}.
See also \href{https://arxiv.org/abs/2312.13335}{arXiv:2312.13335}.

\bibitem{BM}
P.~Bourgade and K.~Mody,
Gaussian fluctuations of the determinant of Wigner matrices,
\emph{Electron. J. Probab.} \textbf{24} (2019), paper no.~96, 1--28.
\href{https://doi.org/10.1214/19-EJP356}{doi:10.1214/19-EJP356}.

\bibitem{CipolloniLandon}
G.~Cipolloni and B.~Landon,
Maximum of the characteristic polynomial of i.i.d.\ matrices,
\emph{Comm. Pure Appl. Math.} \textbf{78} (2025), no.~9, 1703--1782.
\href{https://doi.org/10.1002/cpa.22250}{doi:10.1002/cpa.22250}.

\bibitem{EYY}
L.~Erd\H{o}s, H.-T.~Yau, and J.~Yin,
Rigidity of eigenvalues of generalized Wigner matrices,
\emph{Adv. Math.} \textbf{229} (2012), no.~3, 1435--1515.
\href{https://doi.org/10.1016/j.aim.2011.12.010}{doi:10.1016/j.aim.2011.12.010}.
See also \href{https://arxiv.org/abs/1007.4652}{arXiv:1007.4652}.

\bibitem{ForresterNagaoHonner}
P.~J.~Forrester, T.~Nagao, and G.~Honner,
Correlations for the orthogonal-unitary and symplectic-unitary transitions
at the hard and soft edges,
\emph{Nuclear Phys. B} \textbf{553} (1999), no.~3, 601--643.
\href{https://doi.org/10.1016/S0550-3213(99)00272-2}{doi:10.1016/S0550-3213(99)00272-2}.

\bibitem{Gus}
J.~Gustavsson,
Gaussian fluctuations of eigenvalues in the GUE,
\emph{Ann. Inst. H. Poincar\'e Probab. Statist.} \textbf{41} (2005), no.~2, 151--178.
\href{https://doi.org/10.1016/j.anihpb.2004.04.002}{doi:10.1016/j.anihpb.2004.04.002}.

\bibitem{He}
Y.~He,
Mesoscopic linear statistics of Wigner matrices of mixed symmetry class,
\emph{J. Stat. Phys.} \textbf{175} (2019), no.~5, 932--959.
\href{https://doi.org/10.1007/s10955-019-02266-8}{doi:10.1007/s10955-019-02266-8}.

\bibitem{HuangLandon}
J.~Huang and B.~Landon,
Rigidity and a mesoscopic central limit theorem for Dyson Brownian motion
for general $\beta$ and potentials,
\emph{Probab. Theory Related Fields} \textbf{175} (2019), 209--253.
\href{https://doi.org/10.1007/s00440-018-0889-y}{doi:10.1007/s00440-018-0889-y}.

\bibitem{LandonSosoe}
B.~Landon and P.~Sosoe,
Almost-optimal bulk regularity conditions in the CLT for Wigner matrices,
preprint (2022),
\href{https://arxiv.org/abs/2204.03419}{arXiv:2204.03419}.

\bibitem{LeGall}
J.-F.~Le~Gall,
\emph{Brownian Motion, Martingales, and Stochastic Calculus},
Graduate Texts in Mathematics, vol.~274,
Springer, Cham, 2016.
\href{https://doi.org/10.1007/978-3-319-31089-3}{doi:10.1007/978-3-319-31089-3}.

\bibitem{OR}
S.~O'Rourke,
Gaussian fluctuations of eigenvalues in Wigner random matrices,
\emph{J. Stat. Phys.} \textbf{138} (2010), no.~6, 1045--1066.
\href{https://doi.org/10.1007/s10955-009-9906-y}{doi:10.1007/s10955-009-9906-y}.

\bibitem{PM}
A.~Pandey and M.~L.~Mehta,
Gaussian ensembles of random Hermitian matrices intermediate between
orthogonal and unitary ones,
\emph{Comm. Math. Phys.} \textbf{87} (1983), no.~4, 449--468.
\href{https://doi.org/10.1007/BF01208259}{doi:10.1007/BF01208259}.

\end{thebibliography}
\end{document}